\documentclass{article}
\usepackage{graphicx}

\begin{document}
\title{Stirling's Series Revisited}
\author{Valerio De Angelis}

\date{January 21, 2009}
\maketitle

\section{Introduction}

We present a concise and elementary derivation of the complete asymptotic expansion for the factorial function $n!$,
which we will refer to as Stirling's series.  While 
there have been numerous published proofs of Stirling's series and of its classical dominant term given 
by Stirling's formula
\[\lim_{n\rightarrow \infty} \frac{n! e^n}{n^n\sqrt{2\pi n}} = 1 \]
(see for example   \cite{Nam}, \cite{Mar}, \cite{Pin}, \cite{Kel}), the present treatment produces a new expression for the coefficients. 
In addition, it brings to light the simple relationship between the asymptotic expansions
of $n!$ and $1/n!$ that, even though easily derived from the well-known expansion of $\log \Gamma(z) $ in terms 
of the Bernoulli numbers \cite{WW}, seems to have no simple published proof, and is not mentioned in any of the 
treatments of the Stirling series that we have examined. 

\section{Stirling's formula}
Let \[\phi(z)=2\frac{e^z-1-z}{z^2} = \sum_{j=0}^\infty \frac{2z^j}{(j+2)!}=1+\frac{z}{3}+\cdots,\]
and 
\begin{equation}
\label{1}
f(t)=\frac{1}{\sqrt{\pi}} \int_{-\infty}^\infty e^{-x^2 \phi(xt)} dx, \ \ t>0.
\end{equation}

We show that the integrand in
(\ref{1}) is bounded by an integrable function, uniformly for $t\in[0,1]$.
Since $\phi(u)$ is increasing for $u\geq0$, we have $e^{-x^2 \phi(xt)} \leq e^{-x^2 \phi(0)}=e^{-x^2}$ 
for $x\geq 0$.
We will now show that if $x<0$, then
 $e^{-x^2 \phi(xt)} \leq e^{2x+2}$ for all $t\in [0,1]$. 

Let $\psi(u)= ue^u +u+2-2e^u$. Then $\psi(u)=\sum_{k=3}^\infty \frac{k-2}{k!}u^k$, and 
so $\psi(u)\geq 0$ for $u\geq 0$. Now note that $\psi(-u) = -e^{-u} \psi(u)$, and
\[\frac{\partial}{\partial t} \phi(xt) =\frac{2}{t^3 x^2}\psi(tx).\] So if $x<0$, we find
\[\frac{\partial}{\partial t} \phi(xt) =\frac{2}{t^3 x^2}\psi(-t|x|)=-\frac{2}{t^3 x^2}e^{-t|x|} \psi(t|x|) <0.\] 
Hence $e^{-x^2\phi(xt)}\leq e^{-x^2\phi(x)} =e^{-2e^x} e^{2x+2} \leq e^{2x+2}$ for all
$x<0$ and $t\in[0,1]$. An application of the dominated convergence 
theorem gives us 
\[\lim_{t\rightarrow 0^+} f(t) = f(0)=\sqrt{\pi}/\sqrt{\pi}=1.\] 

Now, \[f(t)=\frac{e^{2/t^2}}{\sqrt{\pi}}\int_{-\infty}^\infty e^{-2e^{xt}/t^2}e^{2x/t} dx, \] 
and the substitution $u=2e^{xt}/t^2$ yields
\[f(t) = \frac{1}{t\sqrt{\pi}}e^{2/t^2}\left(\frac{t^2}{2}\right)^{2/t^2}\int_0^\infty e^{-u} u^{2/t^2-1} du
=\frac{\Gamma\left(2/t^2\right)}{t\sqrt{\pi}}e^{2/t^2}\left(\frac{t^2}{2}\right)^{2/t^2}.\]

Setting $n=2/t^2$ and using $\Gamma(n)=(n-1)!$, we find
\[
f\left(\sqrt{\frac{2}{n}}\right) = \frac{e^n n!}{n^n \sqrt{2\pi n}},
\]
and Stirling's formula is proved.

\section{Stirling's series}

We recall the notion of asymptotic series. If $F(x)$ is a 
(real or complex valued) function of a real variable $x$, we say that  
$\sum_{m=0}^{\infty }c_{m}x^{-m}$
is an 
\emph{asymptotic series} for $F\left( x\right) $ as $x\rightarrow \infty $,
and write 
\[
F(x)\sim \sum_{m=0}^{\infty }c_{m}x^{-m}\ \ \mbox{as}\ \ x\rightarrow \infty ,
\]
if, for each $n\geq 0$, $x^{n+1}\left(F(x)-\sum_{m=0}^{n}c_{m}x^{-m}\right)$
is bounded as $x\rightarrow \infty $. In general, the infinite series $%
\sum_{m=0}^{\infty }c_{m}x^{-m}$ need not converge for any value of $x$.
Nonetheless, the partial sums provide increasingly better approximations to $%
F(x)$ in the limit $x\rightarrow \infty $. Using obvious notation, 
$F(x)\sim G(x)\sum_{m=0}^{\infty }c_{m}x^{-m}$ as $x\rightarrow \infty $ 
means that $\sum_{m=0}^\infty c_m x^{-m}$ is an asymptotic series for $%
F(x)/G(x)$ as $x\rightarrow \infty $.

The integrand in (\ref{1}) is analytic at $t=0$. 
We now obtain the higher-order terms of the asymptotic series for $n!$ (the Stirling's series in the title of this note)
by expanding this analytic function as a Maclaurin series in $t$. In order to bound the remainder, 
we first estimate the size of the derivatives. The estimates of the previous section imply that 
$e^{-x^2\phi(xt)} \leq e^{-2|x|+2}$ for $|x|\geq 2$ and $0\leq t \leq 1$, and from the Maclaurin series for 
$\phi(z)$ we find $\left|\phi^{(k)}(z)\right| \leq |z|^k e^{|z|} $. 
By induction, we find polynomials
$Q_k(x_0,x_1,\ldots,x_k)$ such that 
\[\frac{\partial ^k}{\partial z^k} e^{-x^2\phi(z)} = e^{-x^2\phi(z)} Q_k\left(x,\phi'(z),\phi''(z),\ldots, \phi^{(k)}(z)\right).\]
So we can find a constant $A$ and positive integers $M,m$ (all depending on $k$) such that
\begin{eqnarray}
\left|\frac{\partial ^k}{\partial t^k} e^{-x^2\phi(xt)}\right| &=&  
e^{-x^2\phi(xt)} \left| x^k Q_k\left(x,\phi'(xt),\phi''(xt),\ldots, \phi^{(k)}(xt)\right)\right| \nonumber \\  
 &\leq & A |x|^M e^{m|x|t} e^{-x^2\phi(xt)} \leq A |x|^M e^{m|x|t} e^{-2|x|+2} \nonumber \\
& \leq &  A |x|^M e^{-|x|+2} \label{est}
\end{eqnarray}
for $|x|\geq 2$ and $0\leq t \leq 1/m$.

Write $P_k(x)=x^k Q_k\left(x,\phi'(0),\phi''(0),\ldots,\phi^{(k)}(0)\right)$. Then $P_k(x)$ is a polynomial
in $x$ with the same parity as $k$, and so
\[\int_{-\infty}^\infty e^{-x^2}P_k(x) dx = 0 \ \ \ \mbox{for odd} \ k.\]
Using Leibniz's formula for the $k$th derivative of a product 
and the values $\phi^{(j)}(0)=2/((j+1)(j+2))$, we find that the polynomials $P_k(x)$ satisfy the recurrence relation
\begin{eqnarray*}
P_0(x)&=&1\\
P_{k+1}(x)&=&-2\sum_{j=0}^k \left( \begin{array}{c} k \\ j \end{array}\right)P_j(x)\frac{x^{k-j+3}}{(k-j+2)(k-j+3)}.
\end{eqnarray*}

For each $N\geq 0$, we can write
\begin{eqnarray*}f(t)&=&\frac{1}{\sqrt{\pi}} \sum_{k=0}^{2N+1} \frac{t^k}{k!} \int_{-\infty} ^\infty \left. \frac{\partial ^k}{\partial t^k} e^{-x^2\phi(xt)}\right|_{t=0} dx\\
&+& \frac{1}{\sqrt{\pi}}  \left. \frac{t^{2N+2}}{(2N+2)!} \int_{-\infty} ^\infty\frac{\partial ^{2N+2}}{\partial t^{2N+2}} 
e^{-x^2\phi(xt)}\right|_{t=\xi} dx \\
&=& \frac{1}{\sqrt{\pi}} \sum_{k=0}^{N} \frac{t^{2k}}{(2k)!}\left. \frac{d^{2k}}{dz^{2k}} 
 \int_{-\infty} ^\infty x^{2k} e^{-x^2\phi(z)}dx \right|_{z=0}
+ t^{2N+2} R_N(t)
\end{eqnarray*}
where $0\leq \xi \leq t $. By the estimate (\ref{est}), the integrand in the remainder term is bounded by an integrable function, uniformly for 
small $t$. Hence $R_N(t)$ is bounded as $t\rightarrow 0^+$.

Setting $t=\sqrt{2/n}$, and using the integration formula
\[\int_{-\infty}^{\infty } x^{2k} e^{-cx^{2}}dx=c^{-k-1/2}\frac{\sqrt{\pi }(2k)!}{2^{2k}k!},\]
we find the asymptotic series
\[
\frac{n!e^n}{n^n\sqrt{2\pi n}} \sim \sum_{k=0}^\infty \frac{a_k}{n^k}  \ \mbox{as} \ n\rightarrow \infty 
\]
where
\begin{eqnarray}
a_k&=& \frac{2^k}{(2k)!\sqrt{\pi}}\int_{-\infty}^\infty e^{-x^2}P_{2k}(x)dx   \nonumber \\ 
&=&
\frac{1}{2^k k!}\left[\frac{d^{2k}}{dz^{2k}} 
\left(\frac{z^2}{2(e^z-1-z)}\right)^{k+1/2}\right]_{z=0}.
\label{a2}
\end{eqnarray}

The first few terms of the expansion are
\begin{equation}
\label{3}
 n!\sim \frac{n^{n}\sqrt{2\pi n}}{e^{n}}\left( 1+\frac{1}{12n}+%
\frac{1}{288n^{2}}-\frac{139}{51840n^{3}}-\frac{571}{2488320n^{4}}+\frac{%
163879}{209018880n^{5}}+\cdots \right) .
\end{equation}

\section{Comparing the expansion of $n!$ and $1/n!$}
We now establish the relationship between (\ref{3}) and the asymptotic expansion for $1/n!$, as mentioned in the
introduction. 

Consider the function
\[g(t)=\frac{1}{\sqrt{\pi}} \int_{-\pi/t}^{\pi/t}  e^{-x^2\phi\left(i xt\right)}dx, \ \ t>0. \]
Since $\mbox{Re}(\phi(i\theta))=\left(\frac{\sin(\theta/2)}{\theta/2}\right)^2 $ is monotone decreasing
 on $0\leq \theta \leq \pi$, we have 
\[\mbox{Re}(\phi(i\theta)) \geq \mbox{Re}(\phi(i\pi))=\frac{4}{\pi^2} 
\ \ \mbox{for}\  0\leq \theta \leq \pi.\] 
Hence 
\[ \int_{-\pi/t}^{\pi/t} \left|e^{-x^2\phi\left(i xt\right)}\right|dx
= \int_{-\pi/t}^{\pi/t} e^{-x^2\mbox{\scriptsize Re}(\phi(i xt))}dx 
\leq \int_{-\infty}^{\infty} e^{-4 x^2/\pi^2}dx < \infty ,\]
and so
\[\lim_{t\rightarrow 0^+} g(t)=\frac{1}{\sqrt{\pi}}\int_{-\infty}^\infty e^{-x^2}dx = 1.\]

Now let $t_n=\sqrt{2/n}$. Then
\begin{eqnarray*}
g\left(t_n\right)&=& \frac{1}{t_n\sqrt{\pi}} \int_{-\pi}^{\pi} e^{-\frac{n}{2}\theta^2 \phi(i \theta)}d\theta\\
&=&e^{-n} \frac{1}{t_n\sqrt{\pi}}\int_{-\pi}^{\pi} e^{n e^{i\theta}} e^{-ni\theta}d\theta\\
&=&e^{-n} \frac{1}{t_n\sqrt{\pi}} \sum_{k=0}^\infty \frac{n^k}{k!} \int_{-\pi}^{\pi}e^{i(k-n)\theta} d\theta\\
&=&\frac{n^n\sqrt{2\pi n}}{n!e^n}.
\end{eqnarray*}

\noindent The above calculation is essentially the same as the one used in \cite{Pin}. 

Note that 
\[\left. \frac{d^{2m}}{dt^{2m}}e^{-x^2 \phi(i t x)}\right|_{t=0} =
(-1)^k \left. \frac{d^{2m}}{dt^{2m}}e^{-x^2 \phi(t x)}\right|_{t=0}.\]
Hence, by expanding the integrand of $g(t)$ in a Maclaurin series as before, we find that the 
asymptotic expansion for $1/n!$ is obtained from that of $n!$ by replacing the coefficients
$a_k$ with $(-1)^k a_k$:
\[ \frac{1}{n!}\sim \frac{e^{n}}{n^{n}\sqrt{2\pi n}}\left( 1-\frac{1}{12n}+%
\frac{1}{288n^{2}}+\frac{139}{51840n^{3}}-\frac{571}{2488320n^{4}}-\frac{%
163879}{209018880n^{5}}+\cdots \right) .\]

In particular, the coefficients $a_k$ satisfy the relation
\begin{equation}
\label{conv}
\sum_{k=0}^m(-1)^k a_k a_{m-k} =0, \ \ \ m\geq 1.
\end{equation} 
The convolution relation (\ref{conv}) could also be derived by recalling that the asymptotic expansion
for $\log n!$ in terms of the Bernoulli numbers involves only odd powers of $1/n$ \cite[p.\ 251]{WW} . 
However, the present derivation is based on relatively elementary notions, 
and independently proves the absence of even powers of $1/n$ in the asymptotic series for $\log n!$.

We remark that the functions $f(t)$ and $g(t)$ employed here (once extended as even functions of $t$ 
to a neighborhood of $t=0$) have derivatives of all order at $t=0$, but they are not analytic, because their
Maclaurin series with coefficients $a_k$ and $(-1)^k a_k$ have radius of convergence zero.  
On the other hand, using Cauchy's integral formula to rewrite the $(2k)$th derivative in (\ref{a2}) 
as an integral in the complex plane, we can easily bound the coefficients $a_k$ and prove that 
the exponential generating function for the Stirling's series coefficients $\sum_{k=0}^\infty \frac{a_k}{k!}x^k$
has positive radius of convergence.

\vspace{0.2cm}
\noindent \textbf{Acknowledgments} 

I thank Tewodros Amdeberhan and Victor Moll for helpful discussions, 
Mark Pinsky for suggesting a simplification of the estimates in Section 4, and the referees
for suggestions that resulted in an improved Section 3.

\bibliographystyle{amsplain}

\bibliography{acompat,StirlingSeriesPre}

\noindent Mathematics Department, Xavier University of Louisiana, New Orleans, LA 70125\\
vdeangel@xula.edu
\end{document}